\documentclass[12pt]{article}
\title{Extended finite operator calculus - an example of algebraization of
analysis}
\author{A.K.Kwa\'sniewski*, E.Borak**\\  
\\ *Higher  School of Mathematics and Applied Informatics\\
 PL - 15-021 Bialystok , ul.Kamienna 17,  Poland\\
 e-mail: kwandr@pl
\\** Institute of Computer Science, Bia{\l}ystok University\\
PL-15-887 Bia{\l}ystok, ul.Sosnowa 64, POLAND\\}

\usepackage{amsmath,amsthm,amsfonts,longtable}

\newtheorem{thm}{Theorem}[section]
\newtheorem{prop}{Proposition}[section]
\newtheorem{note}{Note}[section]
\newtheorem{n}{Notation}[section]
\newtheorem{ex}{Example}[section]
\newtheorem{obs}{Observation}[section]
\newtheorem{com}{Comment}[section]
\newtheorem{rem}{Remark}[section]
\newtheorem{defn}{Definition}[section]
\newtheorem{ident}{Identification}[section]
\numberwithin{equation}{section}

\begin{document}
\maketitle
Central European Journal of Mathematics (2005) to
appear

 \vspace{2mm}

\begin{abstract}
``A Calculus of Sequences'' started in 1936 by Ward constitutes
the general scheme for extensions of classical operator calculus
of Rota - Mullin considered by many afterwards and after Ward.
Because of the notation we shall call the Ward`s calculus of
sequences in its afterwards elaborated form - a $\psi $-calculus.

The $\psi $-calculus in parts appears to be almost automatic,
natural extension of classical operator calculus of Rota - Mullin
or equivalently - of umbral calculus of Roman and Rota.

At the same time this calculus is an example of the algebraization
of the analysis - here restricted to the algebra of polynomials.
Many of the results of $\psi $-calculus may be extended to
Markowsky $Q$-umbral calculus where $Q$ stands for a generalized
difference operator, i.e. the one lowering the degree of any
polynomial by one.

This is a review article based on the recent first author
contributions \cite{1}. As the survey article it is supplemented
by the short indicatory glossaries of  notation  and terms  used
by Ward [2], Viskov [7,8] , Markowsky [12], Roman [28-32] on one
side and the Rota-oriented notation on the other side
[9-11,1,3,4,35] (see also \cite{33}).
\end{abstract}
\small{KEY WORDS: extended umbral calculus , Graves-Heisenberg-Weyl algebra\\
MCS (2000) : 05A40 , 81S99}

\vspace{5mm}

{\bf "The modern evolution... has on the whole been marked by a
trend of algebraization. "}
\begin{flushright}
Herman Weyl
\end{flushright}

\section{Introduction}

We shall call the Wards calculus of sequences \cite{2} in its
afterwards last century elaborated form - a $\psi $-calculus
because of the Viskov`s efficient notation \cite{3}-\cite{8}-
adopted from Boas and Buck . The efficiency of the Rota oriented
language and our notation used has been already exemplified by
easy proving of $\psi $-extended counterparts of all
representation independent statements of $\psi $-calculus
\cite{2}. Here these are $\psi $-labelled representations of
Graves-Heisenberg-Weyl (GHW)\cite{3},\cite{1},\cite{16},\cite{17}
algebra of linear operators acting on the algebra $P$ of
polynomials.

As a matter of fact $\psi $-calculus becomes in parts almost
automatic extension of Rota - Mullin calculus \cite{9} or
equivalently - of umbral calculus of Roman and Rota
\cite{9,10,11}. The $\psi $-extension relies on the notion of
$\partial _{\psi} $-shift invariance of operators with $\psi
$-derivatives $\partial _{\psi} $ staying for equivalence classes
representatives of special differential operators lowering degree
of polynomials by one \cite{7,8,12}. Many of the results of $\psi
$-calculus may be extended to Markowsky $Q$-umbral calculus
\cite{12} where $Q$ stands for arbitrary generalized difference
operator, i.e. the one lowering the degree of any polynomial by
one. $Q$-umbral calculus \cite{12} - as we call it - includes also
those generalized difference operators, which are not series in
$\psi $-derivative $\partial_{\psi}  $ whatever an admissible
$\psi $ sequence would be  (for - "admissible" - see next
section).

The survey proposed here  reviews the operator formulation of ``A
Calculus of Sequences'' started in 1936 by Ward \cite{2} with the
indication of the decisive role the $\psi $-representations of
Graves-Heisenberg-Weyl (GHW) algebra account for formulation and
derivation of principal statements of the $\psi $-extension of
finite operator calculus of Rota  and its extensions.

Restating what was said above let us underline  that all
statements of standard finite operator calculus of Rota are valid
also in the case of $\psi $-extension under the almost mnemonic ,
automatic replacement of $\{D,\hat {x}, id \}$ generators of GHW
by their $\psi $-representation correspondents $\{\partial
_{\psi},\hat {x}_{\psi}, id\}$ - see definitions 2.1 and 2.5.
Naturally any specification of admissible $\psi $ - for example
the famous one defining q-calculus - has its own characteristic
properties not pertaining to the standard case of Rota calculus
realization. Nevertheless the overall picture and system of
statements depending only on GHW algebra is the same modulo some
automatic replacements in formulas demonstrated in the sequel. The
large part of that kind of job was already done in \cite{1,3,35}.

The aim of this presentation is to give a general picture ( see:
Section 3) of the algebra of linear operators on polynomial
algebra. The picture that emerges discloses the fact that any
$\psi $-representation of finite operator calculus or equivalently
- any $\psi $-representation of GHW algebra makes up an example of
the algebraization of the analysis with generalized differential
operators \cite{12} acting on the algebra of polynomials.

 We shall delimit all our considerations to the algebra $P$ of polynomials or sometimes to
the algebra of formal series. Therefore the distinction between
difference and differentiation operators disappears. All linear
operators on $P$ are both difference and differentiation operators
if the degree of differentiation or difference operator is
unlimited.

If all this is extended to Markowsky $Q$-umbral calculus \cite{12}
then many of the results of $\psi $-calculus may be extended to
$Q$-umbral calculus \cite{12}. This is achieved under the almost
automatic replacement of $\{D,\, \hat {x}, id\}$ generators of GHW
or their $\psi $-representation $\{\partial _{\psi},\hat
{x}_{\psi},id\}$ by their $Q$-representation correspondents
$\{Q,\hat {x}_{Q} ,id\}$ - see definition 2.5.

The article is supplemented by the short indicatory glossaries of
notation and terms used by Ward \cite{1}, Viskov \cite{7},
\cite{8}, Markowsky \cite{12}, Roman \cite{28}-\cite{31} on one
side and the Rota-oriented \cite{9}-\cite{11} notation on the
other side \cite{3},\cite{4,35,1}.

\section{Primary definitions, notation and general observations}

In the following we shall consider the algebra $P$ of polynomials
$P = ${\bf F}[x] over the field {\bf F}  of characteristic zero. All
operators or functionals studied here are to be understood as {\it linear}
operators on $P$. It shall be easy to see that they are always well defined.

Throughout the note while saying ``polynomial sequence $\left\{
{p_{n}} \right\}_{0}^{\infty}  $'' we mean\\ deg $p_{n}= n$;\;$n
\geq 0$ and we adopt also the convention that deg $p_{n} < 0$  iff
$p_{n} \equiv 0$.

Consider $\Im $ - the family of functions` sequences (in conformity with
Viskov \cite{7},\cite{8},\cite{3} notation ) such that:\\
$\Im = \{\psi;R \supset \left[ {a,b} \right]\;;\;q \in \left[ {a,b}
\right]\;;\;\psi \left( {q} \right):Z \to F\;;\;\psi _{0} \left( {q} \right)
= 1\;;\;\psi _{n} \left( {q} \right) \ne 0;\;\psi _{ - n} \left( {q} \right)
= 0;\;n \in N\}$.\\
We shall call $\psi = \left\{ {\psi _{n} \left( {q} \right)}
\right\}_{n \ge 0} $ ; $\psi _{n} \left( {q} \right) \ne 0$; $n
\ge 0$ and $\psi _{0} \left( {q} \right) = 1$ an admissible
sequence. Let now $n_{\psi}  $ denotes \cite{3,4}
$$
n_{\psi}  \equiv \psi _{n - 1} \left( {q} \right)\psi _{n}^{-1}
\left( {q} \right),n \geq 0.$$

Then  (note that for admissible $\psi $,  $0_{\psi}= 0$)
$$
n_{\psi}! \equiv \psi _{n}^{-1} \left( {q} \right) \equiv
n_{\psi}
\left( {n - 1} \right)_{\psi}
\left( {n - 2} \right)_{\psi}
\left( {n - 3} \right)_{\psi}  ....
2_{\psi}
1_{\psi};
\quad
0_{\psi}!=1
$$

 $n_{\psi} ^{\underline {k}}  = n_{\psi}  \left( {n - 1} \right)_{\psi}
...\left( {n - k + 1} \right)_{\psi}  $,  $\left(
{{\begin{array}{*{20}c}
 {n} \hfill \\
 {k} \hfill \\
\end{array}} } \right)_{\psi}  \equiv \frac{{n_{\psi} ^{\underline {k}}
}}{{k_{\psi}  !}}$ and ${\bf {\exp}}_{\psi}  \{ y\} =
\sum\limits_{k = 0}^{\infty}  {\frac{{y^{k}}}{{k_{\psi}  !}}} $.

\begin{defn}\label{deftwoone}\em{
Let $\psi $ be admissible. Let $\partial _{\psi}  $ be the linear
operator lowering degree of polynomials by one defined according
to $\partial _{\psi } x^{n} = n_{\psi}  x^{n - 1}\;;\; n \geq 0$.
Then $\partial _{\psi}  $ is called the \textit{$\psi
$-derivative}.}
\end{defn}

\begin{rem}\label{remtwoone}
{\em a) For any rational function R the corresponding factorial $
{R\left( {q^{n}} \right)!}$  of the sequence $ R(q^{n})$ is
defined naturally  \cite{3,4,1} as it is defined for $ n_{\psi}$
sequence , i.e. : $R(q^{n})! = R(q^{n})R(q^{n-1})...R(q^{1})$
 The choice $\psi _{n} \left( {q}
\right)$=$\left[ {R\left( {q^{n}} \right)!} \right]^{ - 1}$ and
$R\left( {x} \right) = \frac{{1 - x}}{{1 - q}}$ results in the
well known $q$-factorial\; $n_{q} ! = n_{q} \left( {n - 1}
\right)_{q} !;\quad 1_{q} ! = 0_{q} ! = 1$ while the $\psi
$-derivative $\partial _{\psi}  $ becomes now ($n_{\psi} = n_{q}
$) the Jackson's derivative \cite{25,26,27,2,3} $\partial _{q} $:
\begin{center}
$\left( {\partial _{q} \varphi}  \right)\left( {x} \right) = \frac{{\varphi
\left( {x} \right) - \varphi \left( {qx} \right)}}{{\left( {1-q} \right)x}}$.
\end{center}
b) Note also that if $\psi = \left\{ {\psi _{n} \left( {q}
\right)} \right\}_{n \geq 0} $ and $\varphi = \left\{ {\varphi
_{n} \left( {q} \right)} \right\}_{n \ge 0} $ are two admissible
sequences then [$\partial _{\psi}  $ , $\partial _{\varphi}$]$ =
0$ iff $\psi = \varphi$. Here [,] denotes the commutator of
operators.}
\end{rem}

\begin{defn}\label{deftwotwo}\em{
Let $E^{y}\left( {\partial _{\psi} }  \right) \equiv exp_{\psi} \{
y\partial _{\psi}  \} = \sum\limits_{k = 0}^{\infty}
{\frac{{y^{k}\partial _{\psi}  ^{k}}}{{k_{\psi}!}}} $.\;
$E^{y}\left( {\partial _{\psi} }\right)$ is called the
\textit{generalized translation operator}.}
\end{defn}

\begin{note}\label{notetwoone}
{\em \cite{3,4,1}\\
$E^{a}\left( {\partial _{\psi} }  \right) f(x) \equiv f(x+_{\psi}
a)\;;\; (x +_{\psi} a)^n \equiv E^{a}\left( {\partial _{\psi} }
\right) x^{n}\;;\;E^{a}\left( {\partial _{\psi} }  \right) f =
\sum \limits_{n \geq 0} {\frac{{a^{n}}}{{n_{\psi}  !}}} \partial
_{\psi
}^{n}f$;\\
and in general $(x +_{\psi} a)^n \ne (x +_{\psi} a)^{n-1} (x +_{\psi} a)$.\\
Note also \cite{1} that in general $(1 +_{\psi} (-1))^{2n+1} \ne
0$ ; $n \geq 0$ though $(1 +_{\psi} (-1))^{2n} = 0$; $n \geq 1$.}
\end{note}

\begin{note}\label{notetwotwo}
{\em \cite{1}\\
$\exp_{\psi} \left( {x + _{\psi}  y} \right) \equiv E^x\left(
\partial _{\psi} \right)\exp_{\psi}  \{ y\}$ - while in general
$\exp_{\psi} \{ x + y\} \ne \exp_{\psi}  \{ x\} \exp_{\psi}  \{
y\} $.}
\end{note}
Possible consequent use of the identity $\exp_{\psi} \left( {x +
_{\psi}  y} \right) \equiv \exp_{\psi}  \{ x\} \exp_{\psi}\{ y\} $
is quite encouraging. It leads among others to ``$\psi
$-trigonometry'' either $\psi $-elliptic or $\psi $-hyperbolic via
introducing $\cos_{\psi},\, \sin_{\psi}$ \cite{1},
$\cosh_{\psi}\,,\,\sinh_{\psi}  $ or in general $\psi $-hyperbolic
functions of m-th order $\left\{ {h_{j}^{\left( {\psi} \right)}
\left( {\alpha}  \right)} \right\}_{j \in Z_{m}}  $defined
according to \cite{13}
\[
R \ni \alpha \to \;h_{j}^{(\psi)} \left( {\alpha}  \right) =
\frac{{1}}{{m}}\sum\limits_{k \in Z_{m}}  {\omega ^{ -
kj}\exp_{\psi}\left\{ {\omega ^{k}\alpha}  \right\}\,;\;j \in
Z_{m}},\;\omega = \exp\left\{ {i\frac{{2\pi} }{{m}}} \right\}.
\]
where $1 < m \in N$ and $Z_{m} = \{ 0,1,...,m-1\}$.

\begin{defn}\label{deftwothree}\em{
A polynomial sequence $\left\{ {p_{n}}  \right\}_{o}^{\infty}  $
is of $\psi$ \textit{-binomial type} if it satisfies the
recurrence
\[
E^{y}\left( {\partial _{\psi} }  \right) p_{n} \left( {x} \right) \equiv
p_{n} \left( {x +_{\psi}  y} \right) \equiv \sum\limits_{k \geq 0} {\left(
{{\begin{array}{*{20}c}
 {n} \hfill \\
 {k} \hfill \\
\end{array}} } \right)} _{\psi}  p_{k} \left( {x} \right)p_{n - k} \left(
{y} \right).
\]}
\end{defn}
Polynomial sequences of $\psi $-binomial type \cite{3,4,1} are
known to correspond in one-to-one manner to special generalized
differential operators $Q$, namely to those $Q = Q\left( {\partial
_{\psi} } \right)$ which are $\partial _{\psi}  $-shift invariant
operators \cite{3,4,1}. We shall deal in this note mostly with
this special case,i.e. with $\psi $-umbral calculus. However
before to proceed let us supply a basic information referring to
this general case of $Q$-umbral calculus.

\begin{defn}\label{deftwofour}\em{
Let $P = {\bf F}$[x]. Let $Q$ be a linear map $Q\;:\;P \to P$ such that:\\
$\forall_{p \in P}$ deg $(Qp) = ($deg $p) -1$ (with the convention
deg $p = -1$ means $p = const = 0$). $Q$ is then called a
\textit{generalized difference-tial operator} \cite{12} or
Gel`fond-Leontiev \cite{7} operator.}
\end{defn}

Right from the above definitions we infer that the following holds.

\begin{obs}\label{obstwoone} \em{
Let $Q$ be as in Definition~\ref{deftwofour}. Let $Qx^{n}=
\sum\limits_{k = 1}^{n} {b_{n,k}}  x^{n - k}$ where $b_{n,1} \ne
0$ of course. Without loose of generality take $b_{1,1} = 1$. Then
$\exists\;{\left\{ {q_{k}}  \right\}_{k \geq 2}  \subset {\bf F}}$
and there exists  admissible $\psi $ such that
\begin{equation} \label{eq1}
Q = \partial _{\psi} + \sum\limits_{k \geq 2} {q_{k} \partial
_{\psi} ^{k}}
\end{equation}
if and only if
\begin{equation} \label{eq2}
b_{n,k} = \left( {{\begin{array}{*{20}c}
 {n} \hfill \\
 {k} \hfill \\
\end{array}} } \right)_{\psi}
b_{k,k}; \quad  n \geq k \geq 1;\;b_{n,1} \ne 0;\;b_{1,1} = 1.
\end{equation}
If $\left\{ q_{k} \right\}_k \geq 2$ and an admissible $\psi $
exist then these are unique.}
\end{obs}

\begin{n}\em{
In the case (\ref{eq2}) is true we shall write : $Q = Q\left(
{\partial _{\psi} }  \right)$ because then and only then the
generalized differential operator  Q  is a series in powers of
$\partial _{\psi}  $.}
\end{n}

\begin{rem} \em{
Note that operators of the (\ref{eq1}) form constitute a group
under superposition of formal power series (compare with the
formula (S) in \cite{13}). Of course not all generalized
difference-tial operators satisfy (\ref{eq1}) i.e. are series just
only in corresponding $\psi $-derivative $\partial _{\psi}  $ (see
Proposition~\ref{propthreeone} ). For example \cite{15} let
$Q=\frac{1}{2}D \hat {x} D- \frac{1}{3}D^{3}$.\; Then $Qx^{n} =
\frac{1}{2}n^{2}x^{n - 1} - \frac{1}{3}n^{\underline {3}} x^{n -
3}$ so according to Observation~\ref{obstwoone} $n_{\psi} =
\frac{1}{2}n^{2}$ and there exists no admissible $\psi $ such that
$Q = Q\left( \partial _{\psi}   \right)$.Here  $ \hat {x}$ denotes
the operator of multiplication by  x  while $n^{\underline {k}}$
is a special case of  $n_{\psi}^{\underline {k}}$ for the choice
$n_{\psi} = n$.}
\end{rem}

\begin{obs} \label{obstwotwo} \em{
From theorem 3.1 in \cite{12} we infer that generalized
differential operators give rise to subalgebras $\sum_{Q}$ of
linear maps (plus zero map of course) commuting with a given
generalized difference-tial operator $Q$. The intersection of two
different algebras $\sum_{Q_{1}}$ and $\sum_{Q_{2}}$ is just zero
map added.}
\end{obs}
The importance of the above Observation~\ref{obstwotwo} as well as the
definition below may be further fully appreciated in the context of the
Theorem~\ref{thmtwoone} and the Proposition~\ref{propthreeone} to come.

\begin{defn} \label{deftwofive}\em{
Let $\left\{ {p_{n}}  \right\}_{n \geq 0} $ be the normal
polynomial sequence \cite{12} ,i.e. $p_{0} \left( {x} \right)= 1$
and $p_{n} \left( {0} \right) = 0\;;\;n\geq 1$. Then we call it
the $\psi $-\textit{basic sequence} of the generalized
difference-tial operator $Q$ if in addition $Q\,p_{n} =
n_{\psi}p_{n - 1}$. In parallel we define a linear map $\hat
{x}_{Q} $: $P \to P$ such that $\hat {x}_{Q} p_{n} =\frac{{\left(
{n + 1} \right)}}{{\left( {n + 1} \right)_{\psi} } }p_{n + 1}
;\quad n \geq 0$. We call the operator $\hat {x}_{Q} $ the {\it
dual} to $Q$ operator.}
\end{defn}
When $Q = Q\left( {\partial_{\psi} }  \right) = \partial _{\psi} $
we write for short: $\hat {x}_{Q\left( {\partial _{\psi}} \right)}
\equiv \hat {x}_{\partial _{\psi}}  \equiv \hat {x}_{\psi}$ (see:
Definition \ref{deftwonine}).\\
Of course [$Q,\hat {x}_{Q} $]$ = id$ therefore $\{Q,{\hat {x}}_{Q},
id\}$ provide us with a continuous family of generators of GHW in -
as we call it - $Q$-representation of Graves-Heisenberg-Weyl algebra.\\
In the following we shall restrict to special case of generalized
differential operators $Q$, namely to those $Q = Q\left( {\partial
_{\psi} }  \right)$ which are $\partial _{\psi}  $-shift invariant
operators \cite{3,4,1} (see: Definition \ref{deftwosix}).

At first let us start with appropriate $\psi $-Leibnitz rules for
corresponding $\psi $-derivatives.\\
$\psi${\bf -Leibnitz rules:}

It is easy to see that the following hold for any formal series $f$ and
$g$:\\
for $\partial _{q} $:\;\; $\partial _{q} \left( {f \cdot g} \right) =
\left({\partial _{q} f} \right) \cdot g + \left( {\hat {Q}f} \right) \cdot
\left({\partial _{q} g} \right)$, where $\left( {\hat {Q}f} \right) \left( {x}
\right) = f \left( {qx} \right)$;\\
for $\partial _{R} = R\left( {q\hat {Q}} \right)\partial _{0}
$:\;\; $\partial _{R} (f \cdot g) (z) = R\left( {q\hat {Q}}
\right)\{
(\partial _{0} f) (z)  \cdot  g(z) + f(0) (\partial _{0} g) (z)\}$\\
where - note - $R\left( {q\hat {Q}} \right)x^{n - 1}= n_{R} x^{n - 1}$
; $(n_{\psi}  = n_{R} =n_{R(q)} = R\left( {q^{n}} \right))$ and finally\\
for $\partial _{\psi} =\hat {n}_{\psi}  \partial _{0} $:
\[
\partial _{\psi} (f \cdot g) (z) =\hat {n}_{\psi} \{ (\partial _{o}
f) (z)  \cdot  g(z) + f(0) (\partial _{0} g) (z)\}
\]
where $\hat {n}_{\psi}  x^{n - 1}= n_{\psi}  x^{n - 1}$ ; $n \geq
1$.

\begin{ex} {\em
Let $Q\left( {\partial _{\psi} }  \right)= D \hat {x} D$, where
$\hat {x} f(x)= x f(x)$ and $D = \frac{d}{dx}$.
Then $\psi = \left\{ {\left[ {\left( {n^{2}} \right)!} \right]^{ - 1}}
\right\}_{n \geq 0} $ and $Q = \partial _{\psi}  $.
Let $Q\left( {\partial _{\psi} }  \right) R(q \hat {Q}) \,
\partial _{0} \equiv \partial _{R} $. Then $\psi = \left\{
{\left[ {R\left( {q^{n}} \right)!} \right]^{ - 1}} \right\}_{n \geq 0} $ and
$Q = \partial _{\psi} \equiv \partial_{R} $. Here $R(z)$ is
any formal Laurent series; $\hat {Q}f(x) = f(qx)$
and $n_{\psi} =R(q^{n})$. $\partial _{0} $ is $q = 0$
Jackson derivative which as a matter of fact - being a difference operator
is the differential operator of infinite order at the same time:
\begin{equation}
\label{eq3}
\partial _{0} =
\sum\limits_{n = 1}^{\infty}  {\left( { - 1} \right)^{n + 1}\frac{{x^{n -
1}}}{{n!}}\frac{{d^{n}}}{{dx^{n}}}}.
\end{equation}
Naturally with the choice $\psi _{n} \left( {q} \right)=\left[
{R\left( {q^{n}} \right)!} \right]^{ - 1}$ and $R\left( {x}
\right) = \frac{{1 - x}}{{1 - q}}$ the $\psi $-derivative
$\partial _{\psi}  $ becomes the Jackson's derivative
\cite{25,26,27,2,3} $\partial _{q} $:
\[
\left( {\partial _{q} \varphi}  \right)\left( {x} \right) = \frac{{1-q\hat
{Q}}}{{\left( {1 - q} \right)}}\partial _{0} \varphi \left( {x} \right).
\]
}
\end{ex}

The equivalent to (\ref{eq3}) form of Bernoulli-Taylor expansion
one may find \cite{16} in {\it Acta Eruditorum} from November 1694
under the name ``{\it series univeralissima''}.

(Taylor`s expansion was presented in his ``Methodus incrementorum directa
et inversa'' in 1715 - edited in London).

\begin{defn} \label{deftwosix}\em{
Let us denote by $End(P)$ the algebra of all linear operators
acting on the algebra $P$ of polynomials. Let
$$
\sum\nolimits_{\psi} = \{ T \in End (P);\;\forall \;\alpha \in
{\bf F};\; \left[ {T,E^{\alpha} \left( {\partial _{\psi} }
\right)} \right] = 0 \}.
$$
Then $\sum_{\psi}$ is a commutative subalgebra of $End (P)$ of
${\bf F}$-linear operators. We shall call these operators
$T:\;\partial _{\psi}$\textit{-shift invariant operators}.}
\end{defn}
We are now in a position to define further basic objects of
``$\psi $-umbral calculus'' \cite{3,4,1}.

\begin{defn} \label{deftwoseven}\em{
Let $Q\left( {\partial _{\psi} }  \right):P \to P;$ the linear
operator $Q\left( {\partial _{\psi} }  \right)$ is a $\partial
_{\psi}$-{\it delta operator} iff
\begin{enumerate}
\renewcommand{\labelenumi}{\em \alph{enumi})}
\item $Q\left( {\partial _{\psi} }  \right)$ is $\partial _{\psi}
$ - shift invariant; \item $Q\left( {\partial _{\psi} }
\right)\left( {id} \right) = const \ne 0$  where  id(x)=x.
\end{enumerate}}
\end{defn}

The strictly related notion is that of the $\partial _{\psi}  $-basic
polynomial sequence:

\begin{defn} \label{deftwoeight}\em{
Let $Q\left( {\partial _{\psi} }  \right):P \to P;$ be the
$\partial _{\psi } $-\textit{delta} operator. A polynomial
sequence $\left\{ {p_{n}}  \right\}_{n \ge 0} $; deg
\textit{p}$_{n}$\textit{= n} such that:
\begin{enumerate}
\renewcommand{\labelenumi}{\em \arabic{enumi})}
\item $p_{0} \left( {x} \right) = 1$; \item $p_{n} \left( {0}
\right) = 0$; $n > 0$; \item $Q\left( {\partial _{\psi} }
\right)p_{n} = n_{\psi}  p_{n - 1} $ ,$\partial _{\psi} $-delta
operator $Q\left( {\partial_ \psi} \right)$is called the $\partial
_{\psi}  $-\textit{basic polynomial sequence} of the $\partial
_{\psi}  $-delta operator.
\end{enumerate}}
\end{defn}

\begin{ident} {\em
It is easy to see that the following identification takes place: $\partial
_{\psi}  $-delta operator $Q\left( {\partial _{\psi} }  \right) = \partial
_{\psi}  $-shift invariant generalized differential operator $Q$.
Of course not every generalized differential operator might be considered
to be such.}
\end{ident}

\begin{note}\label{notetwothree}
{\em
Let $\Phi \left( {x;\lambda}  \right) = \sum\limits_{n \geq
0} {\frac{{\lambda ^{n}}}{{n_{\psi}  !}}p_{n} \left( {x} \right)}
$ denotes the $\psi $-exponential generating function of the
$\partial _{\psi}$-basic polynomial sequence $\left\{ {p_{n}}
\right\}_{n \geq 0} $ of the $\partial_{\psi} $-delta operator $Q
\equiv Q\left( {\partial _{\psi}} \right)$ and let $\Phi \left(
{0;\lambda}  \right) = 1$. Then $Q \Phi \left( {x;\lambda} \right)
= \lambda \Phi \left( {x;\lambda} \right)$ and $\Phi $ is the
unique solution of this eigenvalue problem. If in addition (2.2)
is satisfied  then there exists such an admissible sequence
$\varphi $ that $\Phi \left( {x;\lambda} \right)= \exp_{\varphi}
\left\{\lambda x
 \right\} $  (see Example 3.1).}
\end{note}
The notation and naming established by Definitions
\ref{deftwoseven} and \ref{deftwoeight} serve the target to
preserve and to broaden simplicity of Rota`s finite operator
calculus also in its extended ``$\psi $-umbral calculus'' case
\cite{3,4,1}. As a matter of illustration of such notation
efficiency let us quote after \cite{3} the important
Theorem~\ref{thmtwoone} which might be proved using the fact that
$\forall \; Q\left( {\partial _{\psi} }  \right) \quad \exists !$
invertible $S \in \Sigma_{\psi}$ such that $Q\left( {\partial
_{\psi} } \right)= \partial _{\psi} S$. ( For
Theorem~\ref{thmtwoone} see also Theorem 4.3. in \cite{12}, which
holds for operators, introduced by the
Definition~\ref{deftwofive}). Let us define at first what follows.

\begin{defn} \label{deftwonine} {\em (compare with (17) in \cite {8})\\
The \textit{Pincherle $\psi $-derivative} is the linear map ' :
$\Sigma _{\psi} \quad \to  \quad \Sigma _{\psi} $;
\begin{center}
 $T\;$ ' = $T\;\hat {x}_{\psi}  $ \textbf{-} $\hat {x}_{\psi}  T \equiv
 $\textbf{[}$T$, $\hat {x}_{\psi} $\textbf{]}
\end{center}
where the linear map $\hat {x}_{\psi}  :P \to P;$ is defined in
the basis $\left\{ {x^{n}} \right\}_{n \ge 0} $ as follows
\[
\hat {x}_{\psi}  x^{n} = \frac{{\psi _{n + 1} \left( {q}
\right)\left( {n + 1} \right)}}{{\psi _{n} \left( {q}
\right)}}x^{n + 1} = \frac{{\left( {n + 1} \right)}}{{\left( {n +
1} \right)_{\psi} } }x^{n + 1};\quad n \ge 0 .
\]}
\end{defn}

Then the following theorem is true \cite{3}

\begin{thm}\label{thmtwoone}
{\em ($\psi $-Lagrange and $\psi$-Rodrigues formulas \cite{34,11,12,23,3})}\\
Let $\left\{ {p_{n} \left( {x} \right)} \right\}_{n = 0}^{\infty}  $ be
$\partial _{\psi}  $-basic polynomial sequence of the $\partial _{\psi}
$-delta operator $Q\left( {\partial _{\psi} }  \right)$.\\
Let $Q\left( {\partial _{\psi} }  \right) = \partial _{\psi}
S.$ Then for $n>0$:

\begin{enumerate}
\renewcommand{\labelenumi}{\em (\arabic{enumi})}
\item\label{one} $p_{n}(x) = Q\left( {\partial _{\psi} }
\right)$\textbf{'} $S^{-n-1}\;${x}$^{n}$ ;

\item\label{two} $p_{n}(x) = S^{-n}${x}$
^{n} - \frac{{n_{\psi} } }{{n}}$ ($S^{ -
n}\;$)\textbf{'}{x}$^{n-1};$

\item\label{three} $p_{n}(x) = \frac{{n_{\psi} } }{{n}}\hat {x}_{\psi}
S^{ - n}${x}$^{n-1}$;

\item\label{four} $p_{n}(x) = \frac{{n_{\psi} } }{{n}}\hat
{x}_{\psi} (Q\left( {\partial _{\psi} }  \right)$\textbf{'}
)$^{-1} p_{n-1}(x)$  {\em ($\leftarrow$ Rodrigues $\psi $-formula
)}.
\end{enumerate}
\end{thm}

For the proof one uses typical properties of the Pincherle $\psi
$-derivative \cite{3}.Because  $\partial _{\psi}\textbf{'}= id $
we arrive at the simple and crucial observation.

\begin{obs} \label{obstwothree} {\em [3,35]\\
The triples $\{\partial _{\psi},\hat {x}_{\psi}, id \}$ for any
admissible $\psi $-constitute the set of generators of the $\psi
$-labelled representations of Graves-Heisenberg-Weyl (GHW) algebra
\cite{17,18,19,35,1}. Namely, as easily seen [$\partial
_{\psi},\hat {x}_{\psi}  $] $= id$. (compare with
Definition~\ref{deftwofive})}
\end{obs}

\begin{obs} {\em
In view of the Observation~\ref{obstwothree} the general Leibnitz
rule in $\psi $-representation of Graves-Heisenberg-Weyl algebra
may be written (compare with 2.2.2 Proposition in \cite{18}) as
follows
\begin{equation}
\label{eq4}
\partial _{\psi} ^{n}
\quad
\hat {x}_{\psi} ^{m} \quad =
\quad
\sum\limits_{k \ge 0} {\left( {{\begin{array}{*{20}c}
 {n} \hfill \\
 {k} \hfill \\
\end{array}} } \right)\left( {{\begin{array}{*{20}c}
 {m} \hfill \\
 {k} \hfill \\
\end{array}} } \right)k!}\,
\hat {x}_{\psi} ^{m - k}\,
\partial _{\psi} ^{n - k}.
\end{equation}
One derives the above $\psi $-Leibnitz rule from
$\psi$-Heisenberg-Weyl exponential commutation rules exactly the
same way as in $\{D,\hat {x},id \}$ GHW representation - (compare
with 2.2.1 Proposition in \cite{18} ). $\psi $-Heisenberg-Weyl
exponential commutation relations read:
\begin{equation} \label{eq5}
\exp \{t \partial _{\psi} \}\exp \{a \hat{x}_{\psi} \} = \exp \{at
\} \exp \{a \hat{x}_{\psi}  \} \exp \{t
\partial _{\psi} \}.
\end{equation}  }
\end{obs}
To this end let us introduce a pertinent $\psi $-multiplication
$ \ast _{\psi} $ of functions as specified below.

\begin{n} \label{notetwofour} {\textrm{  }\\
$x \ast _{\psi} x^{n} = \hat {x}_{\psi}(x^{n}) = \frac{{\left( {n
+ 1} \right)}}{{\left( {n + 1} \right)_{\psi} } }x^{n + 1};\quad n
\geq 0$ \; hence $x \ast _{\psi} 1 = 1_{\psi}^{-1}x \not
\equiv x$ \\
$x^{n} \ast _{\psi} x = \hat {x}_{\psi} ^{n} ( x ) = \frac{{
1_{\psi}\left({n +1} \right)!}}{{\left( {n + 1} \right)_{\psi}
!}}x^{n + 1} ;\quad n \geq 0$ \;
hence $1 \ast _{\psi} x = x$ therefore\\
$x  \ast _{\psi} \alpha 1 = x  \ast _{\psi} \alpha = \alpha
1_{\psi}^{-1} x$ and $\alpha 1 \ast _{\psi} x = \alpha \ast
_{\psi} x = \alpha x$ and \\
$\forall x,\alpha \in {\bf F};\; f( x)
\ast _{\psi} x^{n} = f(\hat {x}_{\psi})x^{n}$. }
\end{n}

For $k \ne n $ \; $x^{n} \ast _{\psi} x^{k} \ne x^{k} \ast _{\psi}
x^{n}$ as well as $x^{n} \ast _{\psi} x^{k} \ne x^{n+ k}$ - in
general i.e. for arbitrary admissible $\psi $; compare this with
$(x +_{\psi} a)^n \ne (x +_{\psi} a)^{n-1}( x +_{\psi} a)$.\\
In order to facilitate in the future formulation of observations accounted
for on the basis of $\psi $-calculus representation of GHW algebra we shall
use what follows.

\begin{defn}{\em
With Notation \ref{notetwotwo} adopted let us define the \textit{$
\ast _{\psi} $ powers} of x according to

 $x^{n \ast _{\psi} } \equiv $ x $ \ast _{\psi} x^{\left( {n - 1} \right)
\ast _{\psi} } = \hat {x}_{\psi} (x^{\left( {n - 1} \right)\ast
_{\psi }} ) = $ x $ \ast _{\psi} $ x $ \ast _{\psi} $ ... $ \ast
_{\psi} $ x $=\frac{n!}{n_{\psi}  !}x^{n};\quad n \geq 0$.}
\end{defn}
Note that $x^{n\ast _{\psi} }  \ast _{\psi} x^{k\ast _{\psi}
} = \frac{{n!}}{{n_{\psi}  !}} x^{\left( {n + k} \right)\ast _{\psi} }
\ne x^{k\ast _{\psi} }  \ast _{\psi} x^{n\ast _{\psi} }  =
\frac{{k!}}{{k_{\psi}  !}}x^{\left( {n + k} \right)\ast _{\psi} } $ for
$k  \ne n$ and $x^{0\ast _{\psi} }=1$.\\
This noncommutative $\psi $-product $ \ast _{\psi}$ is deviced so as to
ensure the following observations.

\begin{obs} \label{obstwofive} \textrm{ }
\begin{enumerate}
\renewcommand{\labelenumi}{\em (\alph{enumi})}

\item $\partial _{\psi}  x^{n\ast _{\psi} } = n x^{\left( {n - 1}
  \right)\ast _{\psi} } $;\; $n \ge 0$
\item  $\exp_{\psi}[\alpha x] \equiv \exp \{\alpha \hat {x}_{\psi}  \}\bf{ 1}$
\item $\exp [\alpha x]  \ast _{\psi}
(\exp_{\psi} \{\beta \hat {x}_{\psi}  \}\bf {1}) = (\exp_{\psi}
\{[\alpha +\beta ] \hat {x}_{\psi}  \})\bf {1}$
\item $\partial _{\psi} (x^{k} \ast _{\psi} \quad x^{n\ast _{\psi} } )
= (D x^{k}) \ast _{\psi} x^{n \ast _{\psi} } + x^{k} \ast _{\psi}
(\partial _{\psi}  x^{n\ast _{\psi} })$ hence
\item \label{e}
$\partial _{\psi} ( f \ast _{\psi} g) = ( Df) \ast _{\psi} g  + f
\ast _{\psi} (\partial _{\psi} g)$ ; $f,g$  - formal series
\item \label{f}
$f( \hat {x}_{\psi}) g (\hat {x}_{\psi} )$ {\bf 1} $= f(x) \ast
_{\psi} \tilde {g} (x)$ ; $\tilde {g} (x) = g(\hat
{x}_{\psi})${\bf 1}.
\end{enumerate}
\end{obs}
Now the consequences of Leibniz rule (e) for difference-ization of
the product are easily feasible. For example the Poisson $\psi
$-process distribution ${\pi}_{m}(x) = \frac{1}{N(\lambda ,
x)}$p$_m(x);\; \sum\limits_{m \geq 0}\,$p$_{m}(x) = 1$ is
determined by
\begin{equation} \label{eq6}
\textrm{p}_{m}(x) = \frac{{\left( {\lambda x} \right)^{m}}}{{m!}}
\ast _{\psi} \exp_{\psi} [-\lambda x]
\end{equation}
which is the unique solution (up to a constant factor) of the
 $\partial _{\psi}
$-difference equations systems
\begin{equation} \label{eq7}
\partial _{\psi} \textrm{p}_{m}(x) + \lambda  \textrm{p}_{m}(x) = \lambda
\textrm{p}_{m-1}(x) \; m > 0\;;\; \partial _{\psi} \textrm{p}_{0} (x) =
- \lambda  \textrm{p}_{0} (x)
\end{equation}
Naturally $N(\lambda , x) = \exp [\lambda x] \ast _{\psi}
\exp_{\psi} [-\lambda x]$.

As announced - the rules of $\psi $ -product $ \ast _{\psi} $ are
accounted for  on the basis of $\psi $-calculus representation of
GHW algebra. Indeed,it is enough to consult
Observation~\ref{obstwofive} and to introduce $\psi $-Pincherle
derivation $\hat {\partial} _{\psi}  $ of series in powers of the
symbol $\hat {x}_{\psi}  $ as below. Then the correspondence
between generic relative formulas turns out evident.

\begin{obs} {\em
Let $\hat {\partial} _{\psi}  \equiv \frac{{\partial} }{{\partial \hat
{x}_{\psi} } }$ be defined according to $\hat {\partial} _{\psi}
f(\hat {x}_{\psi}) = [\partial _{\psi},f(\hat {x}_{\psi})]$. Then
$\hat {\partial} _{\psi} \hat {x}_{\psi} ^{n} = n \hat {x}_{\psi} ^{n - 1}\;
;\; n \geq 0$ and $\hat {\partial} _{\psi}  \hat {x}_{\psi} ^{n}$ {\bf 1}
$=\partial_{\psi}  x^{n\ast _{\psi} }$ hence {\bf [}$\hat {\partial}_{\psi}
f(\hat {x}_{\psi})${\bf ]1 } $=\partial _{\psi}f(x)$ where $f$ is a
formal series in powers of $\hat {x}_{\psi}  $ or equivalently in $ \ast
_{\psi}$ powers of $x$.}
\end{obs}

As an example of application note how the solution of \ref{eq7} is obtained
from the obvious solution ${\mathbf p}_{m}(\hat {x}_{\psi})$ of the $\hat
{\partial} _{\psi}  $-Pincherle differential equation \ref{eq8} formulated
within G-H-W algebra generated by $\{\partial _{\psi},\hat {x}_{\psi},id \}$
\begin{equation} \label{eq8}
 \hat {\partial} _{\psi} {\mathbf p}_{m} (\hat {x}_{\psi}) + \lambda
 {\mathbf p}_{m}(\hat {x}_{\psi}) = \lambda {\mathbf p}_{m-1}(\hat {x}_{\psi}
 )\, m >0\;;\;\partial _{\psi} {\mathbf p}_{0} (\hat {x}_{\psi}  ) = -
\lambda {\mathbf p}_{0} (\hat {x}_{\psi}. )
\end{equation}
Namely : due to Observation~\ref{obstwofive} (f) p$_{m}(x) = {\mathbf p}
_{m}(\hat {x}_{\psi})${\bf 1}, where
\begin{equation}
{\mathbf p}_{m}(\hat {x}_{\psi}) = \frac{{\left( {\lambda \hat
{x}_{\psi} } \right)^{m}}}{{m!}}\exp_{\psi} [-\lambda \hat
{x}_{\psi}].
\end{equation}

\section{ The general picture  of the algebra $End(P)$
from GHW algebra point of view  }

The general picture from the title above relates to the general
picture of the algebra $End(P)$ of operators on $P$ as in the
following we shall consider the algebra $P$ of polynomials $P$ =
{\bf F}[x] over the field {\bf F}  of characteristic zero.

With series of  Propositions from [1,3,35,21] we shall draw an
over view picture of the situation distinguished by possibility to
develop further umbral calculus in its operator form for {\it any}
polynomial sequences $\left\{ {p_{n}} \right\}_{0}^{\infty} $ [12]
instead of those of traditional binomial type only.

In 1901 it was proved \cite{20} that every linear operator mapping
$P$ into $P$ may be represented as infinite series in operators
$\hat {x}$ and {\it D}. In 1986 the authors of \cite{21} supplied
the explicit expression for such series in most general case of
polynomials in one variable ( for many variables see: \cite{22} ).
Thus according to Proposition 1 from \cite{21} one has:

\begin{prop} \label{propthreeone} {\em
Let $Q$  be a linear operator that reduces by one the degree of each
polynomial.
Let $\{ q_{n} \left( {\hat {x}} \right)\} _{n \ge 0} $ be an arbitrary
sequence of polynomials in the operator $\hat {x}$. Then $\hat {T} =
\sum\limits_{n \ge 0} {q_{n} \left( {\hat {x}} \right)} Q^{n}$
defines a linear operator that maps polynomials into polynomials.
Conversely, if $\hat {T}$ is linear operator that maps polynomials into
polynomials then there exists a unique expansion of the form
\[
\hat {T} = \sum\limits_{n \ge 0} {q_{n} \left( {\hat {x}} \right)} Q^{n}.
\]
}
\end{prop}

It is also a rather matter of an easy exercise to prove the
Proposition 2 from \cite{21}:

\begin{prop} \label{propthreetwo} {\em
  Let $Q$  be a linear operator that reduces by one the
degree of each polynomial. Let $\{ q_{n} \left( {\hat {x}}
\right)\} _{n \ge 0} ${\it be an arbitrary sequence of polynomials
in the operator} $\hat {x}$.  Let a linear operator that maps
polynomials into polynomials be given by
\begin{center}
$\hat {T} = \sum\limits_{n \ge 0} {q_{n} \left( {\hat {x}} \right)} Q^{n}$.
\end{center}
 Let $P\left( {x;\lambda}  \right) = \sum\limits_{n \ge 0}
{q_{n} \left( {x} \right)} \lambda ^{n}$  denotes indicator of
$\hat {T}$. Then there exists a unique formal series $\Phi \left(
{x;\lambda } \right)$;\;$\Phi \left( {0;\lambda} \right) = 1$ such
that:
\begin{center}
 $Q \Phi \left( {x;\lambda}  \right) = \lambda
\Phi \left( {x;\lambda}  \right).$
\end{center}
Then also $P\left( {x;\lambda} \right) = \Phi \left( {x;\lambda}
\right)^{ - 1}\hat {T}\Phi \left({x;\lambda} \right)$. }
\end{prop}

\begin{ex} \em{
Note that \; $\partial _{\psi} \exp_{\psi} \{\lambda x \}= \lambda
\exp_{\psi} \{\lambda x\};\; \exp_{\psi} \left[ {x}
\right]\left| {_{x = 0}}  \right. = 1$. \;\;(*)\\
Hence for indicator of $\hat {T};\; \hat {T} = \sum\limits_{n \ge 0}
{q_{n} \left( {\hat {x}} \right)} \partial _{\psi} ^{n} $ we have:
\begin{center}
 $P\left( {x;\lambda}  \right)= [\exp_{\psi}  \{\lambda \it
{x}\}]^{-1}\hat {T}\exp_{\psi}  \{\lambda {\it x}\}$.\;\;
 (**)
\end{center}
After choosing $\psi _{n} \left( {q} \right)=\left[ {n_{q} !}
\right]^{ - 1}$ we get $\bf {\exp}_{\psi}  \{\it{ x}\} = \exp_{q}
\{x\}$. In this connection note that $\exp_{0} \left( {x} \right)
= \frac{{1}}{{1 - x}}$ and $\exp(x)$ are mutual limit deformations
for $|x|<1$ due to:
\begin{center}
 $\frac{{exp_{0} \left( {z} \right) - 1}}{{z}} = exp_{o} \left( {z} \right)
\Rightarrow exp_{0} \left( {z} \right) = \frac{{1}}{{1 - z}} =
\sum\limits_{k = 0}^{\infty}  {z^{k}};\; |z|<1$ , i.e.
\end{center}
\[
exp\left( {x} \right)\mathrel{\mathop{\kern0pt\longleftarrow}\limits_{1
\leftarrow q}} exp_{q} \left( {x} \right) = \sum\limits_{n = 0}^{\infty}
{\frac{{x^{n}}}{{n_{q}
!}}\;\mathrel{\mathop{\kern0pt\longrightarrow}\limits_{q \to 0}}
\frac{{1}}{{1 - x}}}.
\]
Therefore corresponding specifications of (*) such as $exp_{0}
\left( {\lambda x} \right) = \frac{{1}}{{1 - \lambda x}}$ or
$\exp(\lambda x)$ lead to corresponding specifications of (**) for
divided difference operator $\partial _{0} $ and \it{D} operator
including special cases from \cite{21}.}
\end{ex}

To be complete let us still introduce \cite{3,4} an important
operator $\hat {x}_{Q\left( {\partial _{\psi} }  \right)} $ dual
to $Q\left( {\partial _{\psi} }  \right)$.

\begin{defn} \label{defthreeone} \em{(see Definition~\ref{deftwofive})\\
Let $\left\{ {p_{n}}  \right\}_{n \ge 0} $ be the $\partial
_{\psi} $-basic polynomial sequence of the $\partial _{\psi}
$-delta operator $Q\left( {\partial _{\psi} }  \right)$. A linear
map $\hat {x}_{Q\left( {\partial _{\psi} }  \right)} $: $P \to P$
; $\hat {x}_{Q\left( {\partial _{\psi} }  \right)}\, p_n
=\frac{{\left( {n + 1} \right)}}{{\left( {n + 1} \right)_{\psi} }
}p_{n + 1} ;\quad n \ge 0$ is called the operator {\it dual} to
$Q\left( {\partial _{\psi} }  \right)$.}
\end{defn}

\begin{com} {\em
Dual in the above sense corresponds to adjoint in $\psi$-umbral
calculus language of linear functionals' umbral algebra (compare
with Proposition 1.1.21 in \cite{23} ).}
\end{com}
It is now obvious that the following holds.
\begin{prop}\em{
Let $\{ q_{n} \left( {\hat {x}_{Q\left( {\partial _{\psi} }
\right)}} \right)\} _{n \ge 0} $ be an arbitrary sequence of
polynomials in the operator $\hat {x}_{Q\left( {\partial _{\psi} }
\right)} $. Then $T = \sum\limits_{n \ge 0} {q_{n} \left( {\hat
{x}_{Q\left( {\partial _{\psi} } \right)}}  \right)} Q\left(
{\partial _{\psi} }  \right)^{n}$ defines a linear operator that
maps polynomials into polynomials. Conversely, if $T$ is linear
operator that maps polynomials into polynomials then there exists
a unique expansion of the form
\begin{equation}
\label{eq9}
T = \sum\limits_{n \ge 0} {q_{n} \left( {\hat {x}_{Q\left( {\partial _{\psi
}}  \right)}}  \right)} Q\left( {\partial _{\psi} }  \right)^{n}.
\end{equation}}
\end{prop}

\begin{com} {\em
The pair $Q\left( {\partial _{\psi} }  \right),\; \hat
{x}_{Q\left( {\partial _{\psi} }  \right)} $ of dual operators is
expected to play a role in the description of quantum-like
processes apart from the $q$-case now vastly exploited
\cite{3,4}.}
\end{com}
Naturally the Proposition~\ref{propthreetwo} for $Q\left( {\partial _{\psi} }
\right)$ and $\hat {x}_{Q\left( {\partial _{\psi} }  \right)} $ dual
operators is also valid.\\
{\bf Summing up}: we have the following picture for $End(P)$ -
the algebra of all linear operators acting on the algebra $P$ of
polynomials.
\begin{center}
$Q(P) \equiv \bigcup\limits_{Q} \sum\nolimits_{Q} \subset End(P)$
\end{center}
and of course $Q(P) \ne End(P)$ where the subfamily $Q(P)$ (with zero map)
breaks up into sum of subalgebras $\sum\nolimits_{Q}$ according to
commutativity of these generalized difference-tial operators $Q$ (see
Definition~\ref{deftwofour} and Observation~\ref{obstwotwo}).
Also to each subalgebra $\sum\nolimits_{\psi}$ i.e. to each $Q\left(
{\partial _{\psi} }  \right)$ operator there corresponds its dual operator
$\hat {x}_{Q\left( {\partial _{\psi} }\right)} $
\[
\hat {x}_{Q\left( {\partial _{\psi} }  \right)} \notin
\sum\nolimits_{\psi}
\]
and both $Q\left( {\partial _{\psi} }  \right)$ \& $\hat
{x}_{Q\left( {\partial _{\psi} }  \right)} $ operators are
sufficient to build up the whole algebra $End(P)$ according to
unique representation given by (\ref{eq9}) including the $\partial
_{\psi}  $ and $\hat {x}_{\psi } $ case. Summarising: for any
admissible
$\psi ${\it}  we have the following general statement.\\
{\bf General statement:}
\begin{center}
$End(P) = $[$\{\partial _{\psi}  $,$\hat {x}_{\psi}  $\}] =
[\{$Q\left( {\partial _{\psi} }  \right)$ , $\hat {x}_{Q\left( {\partial
_{\psi} }  \right)} $\}] = [\{$Q$ , $\hat {x}_{Q} $\}]
\end{center}
i.e. the algebra $End(P)$ is generated by any dual pair \{$Q$ ,
$\hat {x}_{Q} $\} including any dual pair \{$Q\left( {\partial _{\psi} }
\right)$ , $\hat {x}_{Q\left( {\partial _{\psi} }  \right)} $\} or
specifically by \{$\partial _{\psi}  $,$\hat {x}_{\psi}  $\} which in turn
is determined by a choice of any admissible sequence $\psi $.

As a matter of fact and in another words: we have bijective
correspondences between different commutation classes of $\partial
_{\psi}  $-shift invariant operators from $End(P)$, different
abelian subalgebras $\sum\nolimits_{\psi}$,\; distinct $\psi
$-representations of GHW algebra, different $\psi
$-representations of the reduced incidence algebra R(L(S)) -
isomorphic to the algebra $\Phi _{\psi}  $ of $\psi $-exponential
formal power series \cite{3} and finally - distinct $\psi $-umbral
calculi \cite{8,12,15,24,34,3,35}. These bijective correspondences
may be naturally extended to encompass also $Q$-umbral
calculi[12,1], $Q$-representations of GHW algebra [1] and abelian
subalgebras $\sum\nolimits_{Q}$.\\
(Recall: R(L(S)) is the reduced incidence algebra of L(S) where\\
L(S)=\{A; A$\subset $S; $|$A$|<\infty $\}; S is countable and
(L(S); $ \subseteq $) is partially ordered set ordered by
inclusion \cite{11,3} ).

This is the way the Rota`s devise has been carried into effect.
The devise {\it ``much is the iteration of the few''} \cite{11} -
much of the properties of literally {\it all} polynomial sequences
- as well as GHW algebra representations - is the application of
few basic principles of the $\psi $-umbral difference operator calculus \cite{3,35,1}.\\
$\psi -$ {\bf Integration Remark :}\\
Recall : $\partial _{o} x^{n} = x^{n - 1}$. $\partial _{o} $ is identical
with divided difference operator. $\partial _{o} $ is identical with
$\partial _{\psi}  $ for $\psi = \left\{ {\psi \left( {q} \right)_{n}}
\right\}_{n \ge 0} $ ; $\psi \left( {q} \right)_{n} = 1$ ; $n \ge 0$ . Let
$\hat {Q}f(x)f(qx)$.\\
Recall also that there corresponds to the ``$\partial _{q} $
difference-ization'' the q-integration \cite{25,26,27} which is a
right inverse operation to ``$q$-difference-ization''\cite{35,1}.
Namely
\begin{equation}
\label{eq10}
F\left( {z} \right): \equiv \left( {\int_{q} \varphi }  \right)\left( {z}
\right): = \left( {1 - q} \right)z\sum\limits_{k = 0}^{\infty}  {\varphi
\left( {q^{k}z} \right)q^{k}}
\end{equation}
i.e.
\begin{multline}
\label{eq11}
F\left( {z} \right) \equiv \left( {\int_{q} \varphi }
\right)\left( {z} \right) = \left( {1 - q} \right)z\left( {\sum\limits_{k =
0}^{\infty}  {q^{k}\hat {Q}^{k}\varphi} }  \right)\left( {z} \right) =\\
=\left( {\left( {1 - q} \right)z\frac{{1}}{{1 - q\hat {Q}}}\varphi}
\right)\left( {z} \right).
\end{multline}
Of course
\begin{equation}
\label{eq12}
\partial _{q} \circ \int_{q} = id
\end{equation}
as
\begin{equation}
\label{eq13}
\frac{{1 - q\hat{Q}}}{{\left( {1 - q} \right)}}\partial _{0} \left(
{\left( {1 - q}\right)\hat {z}\frac{{1}}{{1 - q\hat {Q}}}} \right)=id.
\end{equation}
Naturally (\ref{eq13}) might serve to define a right inverse operation to
``$q$-difference-ization''
\begin{center}
 $\left( {\partial _{q} \varphi}  \right)\left( {x} \right) = \frac{{1 - q\hat
{Q}}}{{\left( {1 - q} \right)}}\partial _{0} \varphi \left( {x} \right)$
\end{center}
and consequently the ``$q$-integration `` as represented by (\ref{eq10}) and
(\ref{eq11}). As it is well known the definite $q$-integral is an numerical
approximation of the definite integral obtained in the $q \to 1$ limit.
Following the $q$-case example we introduce now an $R$-integration
(consult Remark \ref{remtwoone}).
\begin{equation}
\label{eq14}
\int_{R} {x^{n}} = \left( {\hat {x}\frac{{1}}{{R\left( {q\hat {Q}}
\right)}}} \right) x^{n}=
\frac{{1}}{{R\left( {q^{n + 1}} \right)}}x^{n + 1};\quad n \ge 0
\end{equation}
Of course $\partial _{R} \circ \int_{R}= id$ as
\begin{equation}
\label{eq15}
R\left( {q\hat {Q}} \right)\partial _{o}
\left( {\hat {x}\frac{{1}}{{R\left( {q\hat {Q}} \right)}}} \right)=
id .
\end{equation}
Let us then finally introduce the analogous representation for $\partial
_{\psi}  $ difference-ization
\begin{equation}
\label{eq16}
\partial _{\psi}  = \hat {n}_{\psi} \partial _{o};\; \hat {n}_{\psi}
x^{n - 1}= n_{\psi} x^{n - 1};\; n \ge 1.
\end{equation}
Then
\begin{equation}
\label{eq17}
\int_{\psi} {x^{n}} = \left( {\hat {x}\frac{{1}}{{\hat
{n}_{\psi} } }} \right)x^{n}= \frac{{1}}{{\left( {n + 1} \right)_{\psi}
}}x^{n + 1};\;n \ge 0
\end{equation}
and of course
\begin{equation}
\label{eq18}
\partial _{\psi}  \circ \int_{\psi} = id
\end{equation}
{\bf Closing Remark:}

The picture that emerges discloses the fact that any $\psi
$-representation of finite operator calculus or equivalently - any
$\psi $-representation of GHW algebra makes up an example of the
algebraization of the analysis - naturally when constrained to the
algebra of polynomials. We did restricted all our considerations
to the algebra $P$ of polynomials. Therefore the distinction
in-between difference and differentiation operators disappears.
All linear operators on $P$ are both difference and
differentiation operators if the degree of differentiation or
difference operator is unlimited. For example $\frac{{d}}{{dx}} =
\sum\limits_{k \ge 1} {\frac{{d_{k}} }{{k!}}\Delta ^{k}} $ where
$d_{k} = \left[ {\frac{{d}}{{dx}}x^{\underline {k}} } \right]_{x =
0}  = \left( { - 1} \right)^{k - 1}\left( {k - 1} \right)!$ or
$\Delta =\sum\limits_{n \ge 1} {\frac{{\delta _{n}} }{{n!}}}
\frac{{d^{n}}}{{dx^{n}}}$ where $\delta _{n} = \left[ {\Delta
x^{n}} \right]_{x = 0} =1$. Thus the difference and differential
operators and equations are treated on the same footing. For new
applications - due to the first author see [4,1,36-41]. Our goal
here was to deliver the general scheme of "$\psi $-umbral"
algebraization of the analysis of general differential operators
\cite{12}.
 Most of the general features presented here are known to be pertinent to the
{\it Q} representation of finite operator calculus (Viskov,
Markowsky, Roman)  where {\it Q}  is any linear operator lowering
degree of any polynomial by one . So it is most general example of
the algebraization of the analysis for general differential
operators \cite{12}.
\\
\section{Glossary}
In order to facilitate the reader a simultaneous access to quoted
references of classic Masters of umbral calculus - here now follow
short indicatory glossaries of \textbf{notation} used by Ward
\cite{2}, Viskov \cite{7,8}, Markowsky \cite{11}, Roman
\cite{28}-\cite{32} on one side and the \textbf{Rota-oriented
notation} on the other side. See also \cite{33}.
\begin{center}
\begin{longtable}{|c|c|}
          \hline
& \\*
{\bf Ward} & {\bf Rota - oriented} (this note)\\*
& \\*
\endhead
          \hline
& \\*
$\left[n\right];\;\left[n\right]!$ & $n_{\psi};\;n_{\psi}!$\\*
& \\*
basic binomial coefficient $\left[n,r\right]=
\frac{\left[n\right]!}{\left[r\right]!\left[n-r\right]!}$ &
$\psi$-binomial coefficient $\binom{n}{k}_{\psi} \equiv
\frac{n_\psi ^{\underline{k}}}{k_\psi !}$\\*
& \\*
          \hline
& \\*
$D = D_{x}$ - the operator $D$ & $\partial_{\psi}$ - the
$\psi$-derivative\\*
& \\*
$D\,x^n = \left[n\right]\,x^{n-1}$ &
$\partial_{\psi}\,x^n=n_{\psi}\,x^{n-1}$\\*
& \\*
          \hline
&   \\*
$(x+y)^n$ & $(x+_{\psi}y)^n$\\*
&  \\*
$(x+y)^n \equiv \sum\limits_{r=0}^{n}\left[n,r\right]x^{n-r}y^r$  &
$(x+_{\psi} y)^{n}=\sum \limits_{k = 0}^{n}\binom{n}{k}_{\psi}x^ky^{n-k}
$\\*
&  \\*
          \hline

&  \\* basic displacement symbol & generalized shift operator\\* &
\\*$E^t;\;t\in{\bf Z}$ & $E^{y}\left( {\partial _{\psi} }
\right) \equiv exp_{\psi}  \{y\partial _{\psi} \}$;\;$y\in {\bf
F}$\\* &  \\* $E\varphi (x) = \varphi (x+1)$ & $E(\partial
_{\psi})\varphi(x)= \varphi(x+_{\psi}1)$\\* &  \\* $E^t\varphi (x)
= \varphi \left(x+\overline{t}\right)$ & $E^y(\partial _{\psi})x^n
\equiv (x+_{\psi}y)^n$\\* &  \\*
          \hline
& \\*
basic difference operator & $\psi $-difference delta operator\\*
&  \\*
$\Delta = E - id$  & $\Delta _{\psi} = E^y(\partial_{\psi}) - id$\\*
& \\*
$\Delta = \varepsilon (D) - id = \sum\limits_{n=0}^{\infty}
\frac{D^n}{\left[n\right]!} - id$ & \\*
& \\*
          \hline
\end{longtable}
\end{center}
\begin{center}
\begin{longtable}{|c|c|}
          \hline
& \\*
{\bf Roman} & {\bf Rota - oriented} (this note)\\
& \\*
\endhead
          \hline
& \\*
$t;\;tx^n = nx^{n-1}$ & $\partial_{\psi}$ - the $\psi$-derivative\\*
& \\*
& $\partial_{\psi}x^n = n_{\psi}x^{n-1}$\\*
& \\*
$\langle t^{k}|p(x)\rangle = p^{(k)}(0)$ &
$[\partial_{\psi}^kp(x)]|_{x=0}$\\*
& \\*
          \hline
& \\*
evaluation functional & generalized shift operator\\*
& \\*
$\epsilon_{y}(t) = \exp{\{yt\}}$
& $E^y(\partial_{\psi}) = \exp_{\psi}{\{y\partial_{\psi}\}}$\\*
& \\*
$\langle t^k|x^n\rangle = n!\delta_{n,k}$ & \\*
& \\*
$\langle \epsilon_y(t)|p(x)\rangle = p(y)$ &
$[E^y(\partial_{\psi})p_n(x)]|_{x=0} = p_n(y)$\\*
& \\*
$\epsilon_y(t)x^n = \sum\limits_{k \geq 0}\binom{n}{k}
x^ky^{n-k}$ &
$E^y(\partial_{\psi})p_n(x) = \sum \limits_{k \geq 0}\binom{n}{k}_{\psi}
p_k(x)p_{n-k}(y)$\\*
& \\*
          \hline
& \\*
formal derivative & Pincherle derivative \\*
& \\*
$f'(t) \equiv \frac{d}{dt}f(t)$ &
$[Q(\partial_{\psi})]${\bf `}$\equiv
\frac{d}{d\partial_{\psi}}Q(\partial_{\psi})$\\*
& \\*
${\overline f}(t)$ compositional inverse of &
$Q^{-1}(\partial_{\psi})$ compositional inverse of \\*
& \\*
formal power series $f(t)$ & formal power series $Q(\partial_{\psi})$\\*
& \\*
          \hline
& \\*
$\theta_{t};\;\;\theta_{t}x^n = x^{n+1};\;n\geq 0$ &
${\hat x}_{\psi};\;\;{\hat x}_{\psi}x^n =
\frac{n+1}{(n+1)_{\psi}}x^{n+1};\;n\geq 0$ \\*
& \\*
$\theta_{t}t = {\hat x}D$ & ${\hat x}_{\psi}\partial_{\psi} = {\hat x}D =
{\hat N}$\\*
& \\*
          \hline
& \\*
$\sum \limits_{k \geq 0}\frac{s_k(x)}{k_{\psi}!}t^k =$ &
$\sum \limits_{k \geq 0}\frac{s_k(x)}{k_{\psi}!}z^k =$\\*
& \\*
$[g({\overline f}(z))]^{-1}\exp {\{x{\overline f}(t)\}}$ &
$s(q^{-1}(z))\exp_{\psi}{\{xq^{-1}(z)\}}$\\*
& \\*
$\{s_{n}(x)\}_{n \geq 0}$ - Sheffer sequence & $q(t),\,s(t)$ indicators\\*
& \\*
for $(g(t),f(t))$ & of $Q(\partial_{\psi})$ and $S_{\partial_{\psi}}$\\*
& \\*
          \hline
& \\*
$g(t)\,s_{n}(x) = q_{n}(x)$ - sequence &
$s_n(x) = S_{\partial_{\psi}}^{-1}\,q_n(x)$ - $\partial_{\psi}$ - basic\\*
& \\*
associated for $f(t)$ & sequence of $Q(\partial_{\psi})$ \\*
& \\*
          \hline
& \\*
The expansion theorem: & The First Expansion Theorem\\*
& \\*
$h(t) = \sum \limits_{k=0}^\infty \frac{\langle h(t)|p_{k}(x) \rangle}{k!}
f(t)^k$ & $T = \sum \limits_{n\geq 0} \frac{[T\,p_{n}(z)]|_{z=0}}{n_{\psi}}
Q(\partial_{\psi})^n$\\*
& \\*
$p_n(x)$ - sequence associated for $f(t)$ & $\partial_{\psi}$ - basic
polynomial sequence $\{p_n\}_{0}^{\infty}$\\*
& \\*
          \hline
& \\*
$\exp\{y{\overline {f}}(t)\} = \sum \limits_{k=0}^{\infty}\frac{p_{k}(y)}
{k!}t^k$ & $\exp_{\psi}\{xQ^{-1}(x)\} = \sum \limits_{k\geq0}
\frac{p_{k}(y)}{k!}z^k$ \\*
& \\*
          \hline
& \\*
The Sheffer Identity: & The Sheffer ${\psi}$-Binomial Theorem:\\*
& \\*
$s_{n}(x+y) = \sum \limits_{k=0}^{n} \binom{n}{k}p_{n}(y)s_{n-k}(x)$ &
$s_{n}(x+_{\psi} y) = \sum \limits_{k\geq 0}\binom{n}{k}_{\psi}s_{k}(x)
q_{n-k}(y)$\\*
& \\*
          \hline
\end{longtable}
\end{center}
\begin{center}
\begin{longtable}{|c|c|}
          \hline
& \\*
{\bf Viskov} & {\bf Rota - oriented} (this note)\\*
& \\*
\endhead
          \hline
& \\*
$\theta _{\psi }$ - the $\psi$-derivative
& $\partial_{\psi}$ - the $\psi$-derivative\\*
& \\*
$\theta_{\psi}\,x^n = \frac{\psi_{n-1}}{\psi_{n}}x^{n-1}$ &
$\partial_{\psi}\,x^n=n_{\psi}\,x^{n-1}$\\*
& \\*
          \hline
& \\*
$A_p\; (p = \{p_n \}_{0}^{\infty})$  & $Q$\\*
& \\*
$A_p\,p_n = p_{n-1}$ & $Q\,p_n = n_{\psi}p_{n-1}$\\*
& \\*
          \hline
& \\*
$B_p\; (p = \{p_n \}_{0}^{\infty})$  & ${\hat x}_Q$\\*
& \\*
$B_p\,p_n = (n+1)\,p_{n+1}$ & ${\hat x}_Q\,p_n = \frac{n+1}{(n+1)_{\psi}}
p_{n+1}$\\*
& \\*
          \hline
& \\*

$E_{p}^y\; (p = \{p_n \}_{0}^{\infty})$ & $E^{y}\left( \partial
_{\psi} \right) \equiv exp_{\psi}  \{y\partial _{\psi} \}$\\* &
\\* $E_{p}^y p_{n}(x) =\sum \limits_{k=0}^{n}p_{n-k}(x)p_{k}(y)$
& $E^{y}\left( \partial _{\psi} \right)p_{n}(x) = $\\* & \\* &
$=\sum \limits_{k \geq 0}
\binom{n}{k}_{\psi}p_{k}(x)p_{n-k}(y)$\\* & \\*
           \hline
& \\*

$T - \varepsilon_{p}$-operator: & $E^y$ - shift operator:\\* & \\*
$T\,A_{p} = A_{p}\,T$  & $E^y\varphi (x) = \varphi
(x+_{\psi}y)$\\* & \\*
           \hline
& \\*
& $T$ - $\partial_{\psi}$-shift invariant operator:\\*
& \\*
$\forall_{y\in F}\;TE_{p}^y = E_{p}^{y}T$ &
 $\forall_{\alpha \in F}\;[T,E^{\alpha}(\partial_{\psi})]=0$\\*
& \\*
           \hline
& \\*
$Q$ - $\delta_{\psi}$-operator: &
$Q(\partial_{\psi})$ - $\partial_{\psi}$-delta-operator:\\*
& \\*
$Q$ - $\epsilon_{p}$-operator and &
$Q(\partial_{\psi})$ - $\partial_{\psi}$-shift-invariant and\\*
& \\*
$Qx = const \neq 0$ &
$Q(\partial_{\psi})(id) = const \neq 0$\\*
& \\*
          \hline
& \\*
$\{p_n(x),n\geq 0\}$ - $(Q,\psi)$-basic &
$\{p_n\}_{n\geq 0}$ -$\partial_{\psi}$-basic\\*
& \\*
polynomial sequence of the &
polynomial sequence of the\\*
& \\*
$\delta_{\psi}$-operator $Q$ &
$\partial_{\psi}$-delta-operator $Q(\partial_{\psi})$\\*
& \\*
          \hline
& \\*
$\psi$-binomiality property & $\psi$-binomiality property\\*
& \\*
$\Psi_{y}s_n(x) = $ & $E^y(\partial_{\psi})p_{n}(x) = $\\*
& \\*
$=\sum \limits_{m=0}^n \frac{\psi_{n}\psi_{n-m}}{\psi_{n}}s_m(x)
p_{n-m}(y)$ &
$= \sum \limits_{k \geq 0}\binom{n}{k}_{\psi}p_{k}(x)p_{n-k}(y)$\\*
& \\*
          \hline
& \\*
$T = \sum \limits_{n \geq 0}\psi_n[VTp_n(x)]Q^n$ &
$T = \sum \limits_{n\geq 0} \frac{[T p_{n}(z)]|_{z=0}}{n_{\psi}!}
Q(\partial_{\psi})^n$\\*
& \\*
$T\Psi_{y}p(x) =$ & $T p(x+_{\psi} y) =$\\*
& \\*
$\sum \limits_{n \geq 0}\psi_ns_n(y)Q^nSTp(x)$  &
 $\sum \limits_{k \geq 0} \frac{s_{k}(y)}
{k_{\psi}!} Q(\partial_{\psi})^{k}ST p(x)$ \\*
& \\*
          \hline
\end{longtable}
\end{center}

\newpage
\begin{center}
\begin{longtable}{|c|c|}
          \hline
& \\*
{\bf Markowsky} & {\bf Rota - oriented}\\*
& \\*
          \endhead
          \hline
& \\*
$L$ - the differential operator & $Q$\\*
& \\*
$L\,p_n = p_{n-1}$ & $Q\,p_n = n_{\psi}p_{n-1}$\\*
& \\*
          \hline
& \\*
$M$ & ${\hat x}_Q$\\*
& \\*
$M\,p_n = p_{n+1}$  & ${\hat x}_Q\,p_n = \frac{n+1}{(n+1)_{\psi}}
p_{n+1}$\\*
& \\*
          \hline
& \\*
$L_y$ &
$E^{y}\left( Q \right)= \sum \limits_{k \geq 0} \frac{p_k (y)}{k_{\psi}!}
Q^k$\\*
& \\*
$L_y\,p_{n}(x) = $ &
$E^{y}\left( Q \right)\,p_{n}(x) = $\\*
& \\*
$= \sum \limits_{k=0}^{n} \binom{n}{k}p_{k}(x)p_{n-k}(y)$ &
$=\sum \limits_{k \geq 0} \binom{n}{k}_{\psi}\,p_{k}(x)p_{n-k}(y)$\\*
& \\*
           \hline
& \\*
$E^a$ - shift-operator: & $E^y$ - $\partial_{\psi}$-shift operator:\\*
& \\*
$E^a\,f(x) = f(x+a)$ & $E^y\varphi (x) = \varphi (x+_{\psi}y)$\\*
& \\*
           \hline
& \\*
$G$ - shift-invariant operator:
& $T$ - $\partial_{\psi}$-shift invariant operator:\\*
& \\*
$EG=GE$ &
 $\forall_{\alpha \in F}\;[T,E(Q)]=0$\\*
& \\*
           \hline
& \\*
$G$ - delta-operator: &
$L=L(Q)$ - $Q_{\psi}$-delta operator:\\*
& \\*
$G$ - shift-invariant and &
$[L,Q] = 0$ and\\*
& \\*
$Gx = const \neq 0$ & $L(id) = const \neq 0$\\*
& \\*
          \hline
& \\*
$D_L(G)$ & $G\mathbf{'} = [G(Q),{\hat x}_{Q}]$ \\*
& \\*
$L$ - Pincherle derivative of $G$ & $Q$ - Pincherle derivative\\*
& \\*
$D_L(G) = [G,M]$ & \\*
& \\*
          \hline
& \\*
$\{Q_{0},Q_{1},...\}$ - basic family &
$\{p_n\}_{n\geq 0}$ -$\psi$-basic\\*
& \\*
for differential operator $L$ & polynomial sequence of the\\*
& \\*
& generalized difference operator $Q$\\*
& \\*
          \hline
& \\*
binomiality property & $Q$ - $\psi$-binomiality property\\*
& \\*
$P_n(x+y) = $  & $E^y(Q)p_{n}(x) = $\\*
& \\*
$=\sum \limits_{i=0}^n \binom{n}{i}P_i(x)P_{n-i}(y)$ &
$= \sum \limits_{k \geq 0}\binom{n}{k}_{\psi}p_{k}(x)p_{n-k}(y)$\\*
& \\*
          \hline
\end{longtable}
\end{center}
{\bf {\em Acknowledgements:}} The authors thank the Referee for
suggestions , which  have led us to improve the presentation of
the paper. The authors express also their gratitude to Katarzyna
Kwa\'sniewska for preparation the \LaTeX version of this
contribution.\\


\begin{thebibliography}{10}

\bibitem{1}
A. K. Kwa\'sniewski: {\em Bulletin de la Soc. des Sciences et des
Letters de {\L}\'od\'z 52 SERIE Reserchers sur les deformations}
{\bf 36}, 45 (2002).ArXiv: math.CO/0312397

\bibitem{2} M. Ward: {\em Amer. J. Math.\/}
{\bf 58}, 255 (1936).

\bibitem{3}
A. K. Kwa\'sniewski: {\em Rep. Math. Phys.} \textbf{48} (3), 305
(2001) ArXiv: math.CO/0402078  Feb 2004

\bibitem{4}
A. K. Kwa{\'s}niewski: {\em Integral Transforms and Special
Functions}\textbf{2} (4), 333 (2001)

\bibitem{5}
R. P. Boas and Jr. R. C. Buck: {\em Am. Math. Monthly\/} {\bf 63}, 626
(1959).

\bibitem{6}
R. P. Boas and Jr. R. C. Buck: {\em Polynomial Expansions of Analytic
Functions\/}, Springer, Berlin 1964.

\bibitem{7}
O.V. Viskov: {\em Soviet Math. Dokl.\/} {\bf 16}, 1521 (1975).

\bibitem{8}
O.V. Viskov: {\em Soviet Math. Dokl.\/} {\bf 19}, 250 (1978).

\bibitem{9}
G.-C. Rota and R. Mullin: {\em On the foundations of combinatorial
theory, III. Theory of Binomial  Enumeration in "Graph Theory and
Its Applications"\/}, Academic Press, New York 1970.

\bibitem{10}
G. C. Rota, D.Kahaner and A. Odlyzko: {\em J. Math. Anal. Appl.\/} {\bf 42},
684 (1973).

\bibitem{11}
G. C. Rota: {\em Finite Operator Calculus\/}, Academic Press, New York 1975.

\bibitem{12}
G. Markowsky: {\em J. Math. Anal. Appl.\/} {\bf 63}, 145 (1978).

\bibitem{13}
A. K. Kwa\'sniewski: {\em Advances in Applied Clifford Algebras\/} {\bf 9},
41 (1999).

\bibitem{14}
O.V. Viskov: {\em Trudy Matiematicz`eskovo Instituta AN SSSR\/} {\bf 177},
21 (1986).

\bibitem{15}
A. Di Bucchianico and D.Loeb: {\em J. Math. Anal. Appl.\/} {\bf 92}, 1 (1994).

\bibitem{16}
N. Ya. Sonin: {\em Izw. Akad. Nauk\/} {\bf 7}, 337 (1897).

\bibitem{17}
C. Graves: {\em Proc. Royal Irish Academy\/} {\bf 6}, 144 (1853-1857).

\bibitem{18}
P. Feinsilver and R. Schott: {\em Algebraic Structures and Operator
Calculus\/}, Kluwer Academic Publishers, New York 1993.

\bibitem{19}
O.V. Viskov: {\em Integral Transforms and Special Functions\/} {\bf 1},
2 (1997).

\bibitem{20}
S. Pincherle and U. Amaldi: {\em Le operazioni distributive e le loro
applicazioni all`analisi\/},  N. Zanichelli, Bologna 1901.

\bibitem{21}
S. G. Kurbanov and V. M. Maximov: {\em Dokl. Akad. Nauk Uz. SSSR\/} {\bf 4},
8 (1986).

\bibitem{22}
A. Di Bucchianico and D.Loeb: {\em Integral Transforms and Special
Functions\/} {\bf 4}, 49 (1996).

\bibitem{23}
P. Kirschenhofer: {\em Sitzunber. Abt. II  Oster. Ackad. Wiss. Math. Naturw.
Kl.\/} {\bf 188}, 263 (1979).

\bibitem{24}
A. Di Bucchianico and D.Loeb: {\em J. Math. Anal. Appl.\/} {\bf 199}, 39
(1996).

\bibitem{25}
F. H. Jackson: {\em Quart. J. Pure and Appl. Math.\/} {\bf 41}, 193 (1910).

\bibitem{26}
F. H. Jackson: {\em Messenger of  Math.\/} {\bf 47}, 57 (1917).

\bibitem{27}
F. H. Jackson: {\em Quart. J. Math.\/} {\bf 2}, 1 (1951).

\bibitem{28}
S. M. Roman: {\em J. Math. Anal. Appl.\/} {\bf 87}, 58 (1982).

\bibitem{29}
S. M. Roman: {\em J. Math. Anal. Appl.\/} {\bf 89}, 290 (1982).

\bibitem{30}
S. M. Roman: {\em J. Math. Anal. Appl.\/} {\bf 95}, 528 (1983).

\bibitem{31}
S. M. Roman: {\em The umbral calculus\/}, Academic Press, New York 1984.

\bibitem{32}
S. R. Roman: {\em J. Math. Anal. Appl.\/} {\bf 107}, 222 (1985).
\bibitem{33}
A. K. Kwasniewski and E. Gradzka: {\em  Rendiconti del Circolo
Matematico di Palermo Serie II , Suppl.\/} {\bf 69}, 117(2002).
\bibitem{34}
J. F. Steffensen {\em Acta Mathematica\/} {\bf 73}, 333 (1944).
\bibitem{35}
A. K. Kwasniewski: {\em  Integral Transforms and Special
Functions\/} {\bf 14}, 499(2003).



\bibitem{36}
A.K.Kwasniewski    \textit{The logarithmic Fib-binomial formula}
{\em Advan. Stud. Contemp. Math. \/} {\bf 9} No 1 (2004):19-26
ArXiv: math.CO/0406258 13 June 2004

\bibitem{37}
A. K. Kwasniewski \textit{On basic Bernoulli-Ward polynomials}
Bulletin de la Societe des Sciences et des Lettres de Lodz {\bf
54} Serie: Recherches sur les Deformations Vol. 45 (2004) : 5--10
ArXiv: math.CO/0405577    30 May 2004

\bibitem{38}
A. K. Kwasniewski  \textit{$\psi$-Appell  polynomials` solutions
of the -difference calculus  nonhomogeneous equation}  Bulletin de
la Societe des Sciences et des Lettres de Lodz {\bf 54} Serie:
Recherches sur les Deformations Vol. 45  (2004) :  11-15  in print
ArXiv: math.CO/0405578  30 May 2004

\bibitem{39}
A. K. Kwasniewski \textit{On $\psi$-umbral difference
Bernoulli-Taylor formula with Cauchy type remainder} Bulletin de
la Societe des Sciences et des Lettres de Lodz {\bf 54} Serie:
Recherches sur les Deformations Vol. 44  (2004) :21-29 ArXiv:
math.GM/0312401 December 2003

\bibitem{40}
A. K. Kwasniewski \textit{First contact remarks on umbra
difference calculus references streams} , Bull. Soc. Sci. Lett.
Lodz    to appear      ArXiv: math.CO/0403139 v1  8  March 2004
\bibitem{41}
A.K.Kwasniewski  \textit{On extended umbral calculus,
oscillator-like algebras and  Generalized Clifford Algebra},
Advances in Applied Clifford Algebras , {\bf 11} No2
(2001):267-279 , ArXiv: math.QA/0401083 January 2004


\end{thebibliography}
\end{document}